\documentclass[11pt]{article}
\usepackage{amsmath}

\usepackage{amsmath,amsfonts,amssymb,amscd,amsthm}
\usepackage[all]{xy}

\title{Noncommutative identities}
\author{Maxim Kontsevich}

\newcommand{\op}[1]{\operatorname{#1}}
\newcommand{\C}{{\mathbb C}}

\newcommand{\Z}{{\mathbb Z}}

\newcommand{\A}{{\mathcal A}}

\newcommand{\la}{{\langle}}
\newcommand{\ra}{{\rangle}}

\newtheorem{thm}{Theorem}

\newtheorem{conj}{Conjecture}

\begin{document}

\maketitle

 {\raggedleft \it to Don Zagier, on the occasion of his $3\cdot4\cdot 5$ birthday, with love and admiration}
\section{``Characteristic polynomial''}

Let us fix an integer $n\ge 1$ and consider the algebra

$$\A:=\C\la X_1^{\pm 1},\dots ,X_n^{\pm 1}\ra$$
 of noncommutative Laurent polynomials with coefficients in $\C$ in $n$ invertible variables, i.e. the group ring of
 the free group $\op{Free}_n$ on $n$ generators:
$$\A=\C[\op{Free}_n]=\big\{\sum_{g\in \op{Free}_n}  c_g\cdot g\,|\,c_g\in \C,\,\,c_g=0\mbox{ for almost all }g\in \op{Free}_n\big\}\,\,.$$
Define a linear functional $\op{``Tr''}$ on $\A$ by taking the  constant term,
$$\op{```Tr''}:\A\to \C\,,\,\,\,\,\,\,\,\,\,\,\,\,\,\, \sum_g c_g\cdot g\mapsto c_{id}\,\,.$$
This functional vanishes on commutators, like the trace for matrices. By the analogy with the  matrix case, we define the
``characteristic polynomial"\footnote{This is an algebraic analogue of so-called Fuglede-Kadison determinant, see 
\cite{Lueck}.} for any $a\in \A$ as a formal power series in one (central) variable $t$:
$$P_a=P_a(t):=\op{``det"}(1-ta):=\exp\left( -\sum_{k\ge 1}\frac{\op{``Tr"}(a^k)}{k} t^k  \right)=1+\dots\in \C[[t]]\,\,.$$
\begin{thm} For any $a\in \A$ the series $P_a$ is algebraic, i.e.
$$P_a\in \overline{\C(t)}\cap \C[[t]]\subset \overline{\C((t))}\,\,.$$
\end{thm}
Here are few examples: \begin{itemize}
\item the case $n=1$ is elementary, follows easily from the residue formula,
\item for any $n\ge1 $ and
 $$a=X_1+X_1^{- 1}+\dots +X_n+X_n^{-1}$$
one can show that
$$P_a=\frac{\left(\frac{f+1}{2}\right)^n}{\left(  \frac{nf+n-1}{2n-1}\right)^{n-1}}\,,\,\,\,f=f(t):=\sqrt{1-4(2n-1)t^2}=1+\dots\in \Z[[t]]\,.$$
\item if $a=X_1+\dots+X_n+(X_1\dots X_n)^{-1}$ then  the series $P_a$ is an algebraic hypergeometric function.
\end{itemize}
{\it A sketch of the proof}:

Let us assume for simplicity that $a\in \Z[\op{Free}_n]$, the general case is just slightly more complicated.

{\bf Step 1}. For $a=\sum_g c_g\cdot g \in \Z[\op{Free}_n]$ the series $P_a$ also has integer coefficients.
 Indeed,
 it is easy to see that
$$P_a=\prod_{k\ge 1}\prod_{(g_1,\dots,g_k)}(1-c_{g_1}\dots c_{g_k}t^k)$$
where for any $k\ge 1$  we take the product over sequences of elements of $\op{Free}_n$ such that
 $g_1\dots g_k=id$ and the sequence $(g_1,\dots, g_k)$ is strictly smaller than all its cyclic permutations
 for the lexicographic order on sequences associated with some total ordering of $\op{Free}_n$ considered as a countable set.
 Similar argument works if we replace $\op{Free}_n$ by an arbitrary torsion-free group.

 {\bf Step 2}. Consider the following series 
$$F_a=F_a(t):=\sum_{k\ge 1} \op{``Tr"}(a^k)\, t^k\in \Z[[t]].$$
Then $F_a$ is algebraic. This fact  was rediscovered many times, see e.g. \cite{BG}.
 It follows from the theory of noncommutative algebraic series developed
 by N.~Chomsky and M.-P.~Sch\"utzenberger in 1963 (see \cite{Chomsky}). We recommend to the reader to 
 consult chapter 6 in \cite{Stanley},  the algebraicity of the series $F_a$ is the statement of Corollary 6.7.2 in this book.

{\bf Step 3}. Recall the Grothendieck conjecture on algebraicity. It says that for any algebraic vector bundle with flat connection 
over an algebraic variety defined over a number field, all solutions of the corresponding holonomic system of differential equations
  are algebraic if and only if the $p$-curvature vanishes for all sufficiently large primes $p\gg 1$.
 There is a simple sufficient  criterion for such a vanishing. Namely, it is enough to assume that there exists a fundamental system of solutions
 in formal power series at some algebraic point, such that  all Taylor coefficients (in some local algebraic coordinate system) of all solutions have in total  only finitely many primes in denominators.

 The Grothendieck conjecture in its full generality is largely inaccessible  by now. The only two general results are a theorem by N.~Katz on the validity of the 
 Grothendieck conjecture for Gauss-Manin connections, and a theorem of D.~Chudnovsky and G.~Chudnovsky \cite{Chudnovsky}
 in the case of line bundles over algebraic curves (see also Exercise 5 to \textsection 4 of  Chapter  VIII, page 160 in \cite{A}).
   This is exactly our case, by the previous step, and 
because
 $$\frac{d}{dt} P_a=-\frac{F_a}{t} P_a\,.$$ \qed

\section{ Noncommutative integrability, the  case of two variables}

Let now consider the algebra of noncommutative polynomials in two (non-invertible) variables
$$\A=\C\la X,Y\ra\,\,.$$
For any integer $d\ge 1$ we consider the variety $\mathcal{M}_d$ of $\op{GL}_d(\C)$-equivalence classes
 (by conjugation) of $d$-dimensional representations $\rho:\A\to \op{Mat}_{d\times d}(\C)$ of $\A$. More precisely, we are interested
 only in generic pairs of matrices $(\rho(X),\rho(Y))$ and treat variety $\mathcal{M}_d$ birationally. It has dimension $d^2+1$.
   
For generic $\rho$ we consider the ``bi-characteristic polynomial" in two commutative variables 
$$P_\rho=P_\rho(x,y):=\det(1-x\rho(X)-y\rho(Y))=1+\dots \in \C[x,y]\,\,.$$
The  equation $P_\rho(x,y)=0$ of degree $\le d$ defines so-called Vinnikov curve ${\mathcal C}_\rho\subset \C P^2$.
 The number of parameters for the polynomial $P_\rho$ is $\frac{(d+1)(d+2)}{2}-1$, and it is strictly smaller than $\dim  \mathcal{M}_d$ for $d\ge 3$. The missing parameters correspond to the natural line bundle ${\mathcal L}_\rho$ on ${\mathcal C}_\rho$ (well-defined for generic $\rho$) given by 
the kernel of operator $(1-x\rho(X)-y\rho(Y))$ for $(x,y)\in {\mathcal C}_\rho$. Bundle ${\mathcal L}_\rho$ has the same degree as a square root  of the canonical class
 of ${\mathcal C}_\rho$, and defines a point in a torsor over the Jacobian $\op{Jac}({\mathcal C}_\rho)$. For any given generic curve ${\mathcal C}\subset\C P^2$ of degree $d$  line bundles on $\mathcal C$  depend on $\op{genus}({\mathcal C})=\frac{(d-1)(d-2)}{2}$ parameters. Now the dimensions match:
$$\dim \mathcal{M}_d=d^2+1=\frac{(d+1)(d+2)}{2}-1+\frac{(d-1)(d-2)}{2}\,.$$
The conclusion is that $\mathcal{M}_d$ is fibered over the space of planar curves of degree $d$, with the generic fiber being a torsor over an abelian variety. Hence we have  one of simplest examples of an {\it integrable system}.
 Any integrable system has a commutative group of discrete symmetries, i.e. birational automorphisms preserving the structure of the fibration, identical on the base, and acting by shifts on the generic fiber.
 Similarly, one can consider an abelian Lie algebra consisting of  rational vertical vector fields which are infinitesimal generators of shifts on fibers.

Now I want to consider {\it noncommutative} symmetries, i.e. certain universal expressions in free variables $(X,Y)$ which can be specialized
 and make sense for any $d\ge 1$.
An universal  discrete symmetry is an automorphism of $\A$ (or maybe of some completion of $\A$) which preserves the conjugacy class
 of any linear combination $Z(t):=X+tY,\,\,t\in \C$. Indeed, in this case for any $d\ge 1$ and  any representation,  the value of the bi-characteristic
 polynomial $P_\rho$ at any pair of complex numbers $(x,y)\in \C^2$ is preserved, as it can be written as $\det(1-x\rho(Z(y/x)))$.
 Hence the automorphism under consideration is inner on both variables $X$ and $Y$:
$$X\mapsto R\cdot X\cdot R^{-1},\,\,\,Y\mapsto R'\cdot Y\cdot (R')^{-1}\,\,.$$
We are interested in automorphisms of $\A$ only up to inner automorphisms, therefore we may safely assume that $R'=1$.
 Thus, the question is reducing to
 the following one: 

\[\boxed{\begin{aligned}&\mbox{\it  find $R$ such that for any $t\in \C$ there exists $R_t$ such that  }\\
 &\,\,\,\,\,\,\,\,\,\,\,\,\,\,\,\,\,\,\,\,\,\,\,R\cdot X\cdot R^{-1}+tY=R_t\cdot(X+tY)\cdot R_t^{-1}\,\,
\end{aligned}}\]

First, let us make  calculations on the Lie level. Denote by $\mathfrak{g}$ the Lie algebra of derivations $\delta$ of $\A$ of the form
  $$\delta(X)=[D,X] \mbox{ for some } X\in \A,\,\,\,\,\delta(Y)=0$$ 
and such that for any $t\in \C$ there exists $D_t\in \A$ such that
$$\delta(X+tY)=[D_t,X+tY]\,\Longleftrightarrow \, [D,X]=[D_t,X+tY]\,.$$
 It is easy to classify such derivations, and one can check that the following elements  form a linear basis of $\mathfrak{g}$:
$$\delta_{n,m}(X)=[c_{n,m},X],\,\,\delta_{n,m}(Y)=0,\,\,\,n\ge 0,\,m\ge 1$$
where for any $n,m\ge 0$ we define
$$c_{n,m}:=\sum_{\frac{(n+m)!}{n!m!} \,\op{shuffles}\,\,w} w\,,$$
i.e. the sum of all words in $X,Y$ containing $n$  letters  $X$ and $m$  letters $Y$.
 Elements $D_t\in \A$ corresponding to the derivation $\delta_{n,m}$ are given by
$$D_t=\sum_{ 0\le k\le n}c_{n-k,m+k}\,t^k\,\,.$$
A direct calculation shows that $\mathfrak{g}$ is an abelian Lie algebra.

Let us go now the completions of algebra $\A$, and of Lie algebra $\mathfrak{g}$: 
$$\widehat{\A}:=\C\la\la X,Y\ra\ra,\,\,\,\,\widehat{\mathfrak{g}}:=\prod_{n\ge 0,m\ge 1} \C\cdot \delta_{n,m}\,\,.$$
Then the action of $\widehat{\mathfrak{g}}$ on $\widehat{\A}$ exponentiates a continuous group action
$$\widehat{\mathfrak{g}}\stackrel{\op{exp}}{\simeq} \widehat{G}\subset \op{Aut}(\widehat{\A})\,\,.$$
For any $\delta\in \widehat{\mathfrak{g}}$ the corresponding one-parameter group of automorphisms acts by
$$\op{exp}(\tau\cdot \delta):  X\mapsto R(\tau)\cdot X\cdot R(\tau)^{-1},\,\,\,\,Y\mapsto Y\,\,\,\,\,\,\,\,\,\forall\,\tau\in \C$$
for certain invertible element $R(\tau)\in \widehat{\A}^\times$.
An easy calculation shows that $R(\tau)$ is the unique solution of the differential equation
$$\frac{d}{d\tau} R(\tau)=\delta(R(\tau))+R(\tau)\cdot D,\,\,\,\,\,\,\,R(0)=1$$
where $D\in \widehat{\A}$ is such that $\delta(X)=[D,X]$. The value   $R(\tau)_{|\tau=1}$ gives $\exp(\delta)$. 

Now we can start to look for a class of elements $\delta \in \widehat{\mathfrak{g}}$
 such that the corresponding automorphism $\op{exp}(\delta)$ is sufficiently nice, e.g. if it makes some sense for $\A$ without passing to the completion. 

 Let us encode a generic element $\delta$ as before by the corresponding generating series in {\it commutative}
 variables $x,y$:
$$\delta=\sum_{n,m} f_{n,m} \delta_{n,m}\in \widehat{\mathfrak{g}}\,\,\, \leadsto\,\,\, \widetilde{\delta}:=\sum_{n,m} f_{n,m} x^n y^m\in \C[[x,y]]\,\,.$$
I suggest the following Ansatz:
$$
\boxed{\mbox{\it  $\widetilde{\delta}$ is the logarithm of a rational function in $x,y$.}}
$$
Hypothetically, for such $\delta$ the corresponding automorphism $\exp(\delta)$  of $\widehat{\A}$ can be extended to certain ``algebraic extension'' of $\A$.
A good indication is
\begin{thm}For any $P=P(x,y)=1+\dots\in \C[x,y]$ expand 
$$\log(P)=\sum_{n,m} f_{n,m}x^n y^m\in \C[[x,y]]\,\,.$$
Then the series 
$$\exp\left(\sum_{n,m} \frac{(n+m)!}{n! m!} f_{n,m} x^n y^m\right)$$
is algebraic.
\end{thm}
This result is  elementary, and I leave it as an exercise to the reader. (Hint: use the residue formula twice.)
It implies that the image of $R$ under the abelianization morphism
$$\C\la\la X,Y \ra\ra \twoheadrightarrow \C[[x,y]],\,\,\,X\mapsto x,\,Y\mapsto y$$
is algebraic.

 {\bf Example}. Consider the case 
$$\widetilde{\delta}=\log(1-xy)=-\sum_{k\ge 1} \frac{x^k y^k}{k}\,\,.$$
Then one can show that
$$R=1-YX-C\in \widehat{\A}^\times$$
where  $C$   
 is the unique solution of the equation
$$C=X\cdot(1-C)^{-1}\cdot Y\,\,.$$
It can be written 
$$C=XY+XXYY+XXYXYY+XXXYYY+\dots=()+(())+(()())+((()))+\dots$$
as the sum of all irreducible bracketings  if we replace $X$ by $($ and $Y$ by $)$. 

The equation for $C$ is equivalent to the generic ``quadratic equation''
$$T^2+AT+B=0$$
by the  substitutions
$$A:=X^{-1},\,\,\,\,B:=-X^{-1}Y,\,\,\,\,\,T:=X^{-1}C\,\,.$$
The invertible elements $R_t\in \widehat{\A}^\times,\,\,t\in \C$ are given by
$$R_t:=R\cdot(1-t T^2),\,\,\,\,\,R\cdot X\cdot R^{-1}+tY=R_t\cdot (X+tY)\cdot R_t^{-1}\,\,.$$

I'll finish with another example of an integrable system.
Few years ago together  with S.~Duzhin we discovered numerically that the rational map
$$S_{-1}:(X,Y)\mapsto (XYX^{-1},\,(1+Y^{-1})\,X^{-1})$$
should be a discrete symmetry of an integrable system, where $X,Y$ are two $d\times d$ matrices for $d\ge 1$.
Recently O.~Efimovskaya and Th.~Wolf found an explanation (see \cite{EW}). Namely, their results suggest that the conjugacy class
 of the Lax operator which is the matrix $L(t)$ of size $2d\times 2d$, defined as
$$L(t):=\begin{pmatrix}
Y^{-1}+X & tY+Y^{-1}X^{-1}+X^{-1}+1 \\
Y^{-1}+\frac{1}{t}X & Y+Y^{-1}X^{-1}+X^{-1}+\frac{1}{t}
\end{pmatrix}$$ 
does not change under the discrete symmetry $S_{-1}$  as above, for any $t\in \C$. Indeed, one can check directly that
$S_{-1}(L(t))=V(t) L(t) V(t)^{-1}$ where\footnote{I am grateful to A.Odesskii for  help in finding the  matrix $V(t)$.} 
$$V(t):=\begin{pmatrix}
X(1+X+Y)^{-1} X (1+Y^{-1})^{-1} & t X(1+X+Y)^{-1} Y\\
X(1+Y^{-1})^{-1} X (1+X+Y)^{-1} & XY(1+X+Y)^{-1}
\end{pmatrix}\,\,.$$

\section {Noncommutative integrability for many variables}

Let $M=(M_{ij})_{1\le i,j\le 3}$ be a matrix whose entries are $9=3\times 3$ free independent noncommutative variables.
Let us consider 3 ``birational involutions" 
$$\begin{aligned}
&I_1: M \mapsto  M^{-1}\\
&I_2: M \mapsto  M^t\\
&I_3: M_{ij}  \mapsto (M_{ij})^{-1} \,\,\,\,\,\,\forall i,j\,\,.
\end{aligned}$$
The composition $I_1\circ I_2\circ I_3$ commutes with the multiplication on the left and on the right by diagonal
 $3\times 3 $ matrices. We can factorize it by the action of $\op{Diag}_{\op{left}}\times \op{Diag}_{\op{right}}$ and get only 
$4$ independent variables, setting e.g. $M_{ij}=1$ for $i=3$ and/or $j=3$.

\begin{conj} The transformation $(I_1\circ I_2\circ I_3)^3$ is equal to the identity modulo $\op{Diag}_{\op{left}}\times \op{Diag}_{\op{right}}$-action. In other words, there exists two diagonal $3\times 3$ matrices $D_L(M),D_R(M)$
whose entries are noncommutative rational functions in $9$ variables $(M_{ij})$, such that
$$(I_1\circ I_2\circ I_3)^3(M)=D_L(M)\cdot M\cdot D_R(M)\,\,.$$
\end{conj}

This is a very degenerate case of integrability. 
Similarly, for $4\times 4$ matrices
 the transformation $I_1\circ I_2\circ I_3$ should give a genuinely nontrivial integrable system. In the simplest case
 when the entries of this matrix are scalars, the Zariski closure of the generic orbit (modulo the left and the right diagonal actions) is a union of two elliptic curves.

Finally, I'll present a series of hypothetical discrete symmetries of  integrable systems written as recursions.
 Fix an odd integer $k\ge 3$ and consider sequences
$(U_n)_{n\in \Z}$ (of, say, $d\times d$ matrices), satisfying
$$\begin{aligned}
U_n=\, &U_{n-k}^{-1}\, (1+U_{n-1}\, U_{n-k+1}) &\mbox{ for }n\in &\,2\,\Z \\
 U_n= \,   &(1+U_{n-k+1}\,U_{n-1})\, U_{n-k}^{-1}  &\mbox{ for }n\in &\,2\,\Z+1\,.
\end{aligned}$$
Then the map $(U_1,\dots,U_k)\mapsto (U_3,\dots,U_{k+2})$ is integrable.
\section{Noncommutative Laurent phenomenon}

In the previous example one observes also the noncommutative Laurent phenomenon:
$$\forall \,n\in \Z\,\,\,\,\,\,\,\,\,U_n\in \Z\la U_1^{\pm 1},\dots, U_k^{\pm 1}\ra\,.$$
Also with S.~Duzhin we discovered that the noncommutative birational map
$$S_{\,l}:(X,Y)\mapsto (XYX^{-1}, (1+Y^{\,l})\,X^{-1})$$
for $l\ge 1$ satisfies the same Laurent properties, i.e. both components of 2-dimensional vector obtained by an arbitrary
number of iterations, belong to the ring $\Z\la X^{\pm 1},Y^{\pm 1}\ra$.
 The case $l=1$ is easy, and the case $l=2$ was studied by A.~Usnich (unpublished) and by Ph.~Di Francesco and R.~Kedem, see \cite{DF}.
The Laurent property has now three different proofs for the case $l\ge 3$ when the dynamics is {\it non-integrable}:
\begin{itemize}
\item by A.~Usnich using triangulated categories, see \cite{Usnich},
\item an elementary algebraic  proof by A.~Berenstein and V.~Retakh, see \cite{Retakh},
\item a new combinatorial proof of Kyungyong Lee, which also  shows  that all the coefficients of noncommutative Laurent polynomials
obtained by iterations, belong to $\{0,1\}\subset \Z$, see \cite{Lee}. 
\end{itemize}

Finally, recently A.~Berenstein and V.~Retakh found a large class of noncommutative mutations related with triangulated surfaces, and proved the noncommutative Laurent property for them.

 IHES, 35 route de Chartres, F-91440, France

{maxim@ihes.fr}\\

\end{document}